\documentclass[12pt]{amsart}

\input xy
\xyoption{all}
\usepackage{latexsym,amsfonts,amsmath,amssymb,amsthm,rotating,txfonts,mathrsfs}
\usepackage{epic}
\usepackage{t1enc}
\usepackage[latin2]{inputenc}
\usepackage{diagrams}
\newtheorem{thm}{Theorem}
\newtheorem{prop}{Proposition}
\newtheorem{conjecture}{Conjecture}
\newtheorem{defi}{Definition}

\newtheorem{cor}{Corollary}
\newtheorem{exa}{Example}

\newcommand{\jetreg}[2]{J_k^{\mathrm{reg}}({#1},{#2})}

\newcommand{\ann}{\mathrm{Ann}}
\newcommand{\proj}{\mathrm{Proj}}
\newcommand{\tp}{\mathrm{Tp}}
\newcommand{\thom}{\mathrm{Thom}}

\newcommand{\mdeg}[1]{\mathrm{mdeg}[#1]}

\newcommand{\sires}{\res_{\mathbf{z}=\infty}}
\newcommand{\symk}{\Sym^{\le k} \C^n}
\newcommand{\grass}{\mathrm{Grass}}
\newcommand{\flag}{\mathrm{Flag}}
\newcommand{\reg}{\mathrm{reg}}

\newcommand{\C}{\mathbb{C}}
\newcommand{\Q}{\mathbb{Q}}
\newcommand{\PP}{\mathbb{P}}

\newcommand{\calx}{\mathcal{X}}
\newcommand{\calo}{\mathcal{O}}
\newcommand{\cotx}{T_X^*}
\newcommand{\res}{\operatornamewithlimits{Res}}
\newcommand{\dbz}{\,d\mathbf{z}}
\newcommand{\bz}{\mathbf{z}}
\newcommand{\bi}{\mathbf{i}}
\newcommand{\Hom}{\mathrm{Hom}}
\newcommand{\Sym}{\mathrm{Sym}}
\newcommand{\hm}[2]{\mathrm{Hom}(\C^{#1},\C^{#2})}

\title{Moduli of map germs, Thom polynomials and the  Green-Griffiths conjecture}

\author{Gergely B\'erczi -- Oxford}

\begin{document}

\begin{abstract}
This survey paper---based on my IMPANGA lectures given in the Banach Center, Warsaw in January 2011---studies the moduli of holomorphic map germs from the complex line into complex compact manifolds with applications in global singularity theory and the theory of hyperbolic algebraic varieties.\\


\noindent MSC 55N91, 14E15, 32Q45
\end{abstract}

\maketitle

\section{Introduction}

Let $J_k(n,m)$ denote the complex vector space of $k$-jets of map germs from $\C^n$ to $\C^m$ mapping the origin to the origin. The open dense subset $J_k^\reg(n,m)$ consists of jets with regular linear part. $J_k^\reg(1,1)$ is a group under composition of jets, and it acts via reparametrisation on $J_k(1,n)$.  

The dimension of the complex vector space $J_k(1,1)$ is $k$, and with a natural choice of basis $J_k^\reg(1,1)$ can be identified with the following linear subgroup of $GL(k)$:
\[\mathbf{G}_k=\left\{\left(
\begin{array}{ccccc}
\alpha_1 & \alpha_2   & \alpha_3          & \ldots & \alpha_k\\
0        & \alpha_1^2 & 2\alpha_1\alpha_2 & \ldots & 2\alpha_1\alpha_{k-1}+\ldots \\
0        & 0          & \alpha_1^3        & \ldots & 3\alpha_1^2\alpha_ {k-2}+ \ldots \\
0        & 0          & 0                 & \ldots & \cdot \\
\cdot    & \cdot   & \cdot    & \ldots & \alpha_1^k
\end{array}
 \right): \alpha_1 \in \C^* , \alpha_i \in \C  \right\}
,\]
where the polynomial in the $(i,j)$ entry is
\[p_{i,j}(\bar{\alpha})=\sum_{a_1+a_2+\ldots +a_i=j}\alpha_{a_1}\alpha_{a_2} \ldots \alpha_{a_i}.\]

This paper is an exploration of this subgroup of $GL_k$ and the non-reductive quotient $J^k(1,n)/J_k(1,1)$, which is roughly speaking the moduli of $k$-jets of curves in $\C^n$. Principles and ideas of classical reductive geometric invariant theory of Mumford do not apply in this situation, for more details about the background see \cite{dk,bk}.

We illustrate the importance of this moduli space for two classical problems. 

The first problem goes back to Ren\'e Thom \cite{thom} and his study of degreneracy loci of holomorphic maps between manifolds. Consider a holomorphic map $f:N\to M$ between two complex manifolds,
of dimensions $n\leq m$. For a singularity class $O \subset J_k(n,m)$ we can define the set
\[ Z_O[f] = \{p\in N;\; f_p\in O\}, \] 
that is the set of points where the germ $f_p$ belongs to $O$. Then, under some additional technical
assumptions, $Z_O[f]$ is an analytic subvariety of $N$.  The
computation of the Poincar\'e dual class $\alpha_O[f]\in H^*(N,\mathbb{Z})$ of
this subvariety is one of the fundamental problems of global singularity
theory. It turns out that these classes---the Thom polynomials of singularities---are certain equivariant intersection numbers on the moduli space $J_k(1,n)/J_k^\reg(1,1)$.  See \cite{bsz, damon, FR2, FR3, rf, rimanyi, ronga, gaffney, kazarian, prag, prag2, thom}  for details on Thom polynomial computations.

The second problem is an old conjecture of Green and Griffiths about holomorphic curves in smooth projective varieties. Their conjecture, from 1979, says that any projective variety $X$ of general type contains a proper subvariety $Y\subsetneq X$ such that any entire holomorphic curve $f:\C \to X$ sits in $Y$, that is $f(\C) \subset Y$. The strategy of Green-Griffiths \cite{gg}, Demailly \cite{dem} and Siu \cite{siu1,siu2,siu3}, and the recent work of Diverio, Merker and Rousseau \cite{dmr} leads us to prove the positivity of an intersection number on the Demailly bundle, whose fibres are canonically isomorphic to $J_k(1,n)/J_k^\reg(1,1)$. For details of this approach see \cite{ahlfors, bloch, dem, dr, dmr, berczi, merker3, dem3, siu2, siu3}.  

This survey paper is an extended version of my IMPANGA lectures given in the Banach Center, Warsaw in January 2011. I would like to thank to Piotr Pragacz for the warm welcome there. 

Most results presented here have already been published in the papers \cite{bsz, berczi, bk}. The only exception is the formula for the Euler characteristic of Demailly jet bundles in \S8 Appendix and the relation to the curvilinear Hilbert scheme in the last section.

\section{Equivariant cohomology}

It is well-known that any group action on a topological space carries topological information about the space. 

Let $G$ be a topological group. A principal $G$-bundle is a map $E \to B$, which is locally a projection $U \times G \to U$. One of the main fundamental principles in topology is to find universal objects such that all objects in a given category can be "pulled-back" from this. Here a universal principal $G$-bundle is a bundle $\pi:EG \to BG$ such that every principal $G$-bundle $E \to B$ is a pull-back via a map $B\to BG$, which is unique up to homotopy.  
$EG$ is contractible. In fact, if $P$ is a contractible space with a free $G$-action then $P \to P/G$ is a universal principle $G$-bundle. 
   
\begin{thm}
$EG$ exists for all topological group $G$, and unique up to equivariant homotopy.  
\end{thm}

\begin{exa} $B\C^*=\PP^{\infty}(\C)$, and $\C^\infty \to \PP^{\infty}(\C)$ is a universal principle $\C^*$ bundle. 

Similarly, \[BGL_n=\Hom(\C^n,\C^\infty)/GL_n=\mathrm{Gr}(n,\infty),\]
 and $EGL(n)=\Hom(\C^n,\C^\infty) \to \mathrm{Gr}(n,\infty)$ is the universal principle $GL(n)$-bundle. From this we can construct   
\[EGL_n \times_{GL_n} \C^n \to BGL_n,\] 
which is a universal vector bundle, namely any vector bundle of rank $n$ can be pulled back from this.
\end{exa}

The next step is to define equivariant cohomology.
Let $X$ be a $G$-space, i.e. a topological space with a $G$-action. If the action is free, then $G$-equivariant cohomology is the ordinary cohomology of the quotient $H^*(X/G)$. For non-free actions the quotient $X/G$ is not well-behaved and $H^*(X/G)$ does not carries enough information. We need to "resolve" the action by replacing $X$ with $X \times EG$. This has a free (diagonal) $G$-action, and define  
\[H^*_G(X)=H^*(EG\times_G X)\]
 
\begin{exa}
$H_G^*(pt)=H^*(BG)=\C[\mathfrak{h}]^W$, where $\mathfrak{h}=Lie T$ is the Cartan algebra acted on by the Weil group $W$. For example $H^*_{GL_n}(pt)=S^W=\C[x_1,\ldots, x_n]^{S_n}$, the algebra of symmetric polynomials.  
\end{exa}

\noindent \textbf{Properties of equivariant cohomology:} 
\begin{enumerate}
\item $f:X \to Y$ $G$-map induces $H(f):H_G(Y) \to H_G(X)$  
\item $h:G\to H$ homomorphism, then $EH$ can serve as $EG$ and we have a projection $EH \times_G X \to EH \times_H X$ which induces $H(h):H_H(X) \to H_G(X)$
\item $H_G^*(pt)=H^*(BG)=\C[\mathfrak{h}]^W$, and $H_G^*(X)$ is a $H_G^*(pt)$-module. For example $H^*_{GL_n}(pt)=S^W=\C[x_1,\ldots, x_n]^{S_n}$. 
\end{enumerate}

\begin{prop}{\textbf{Induction, Restriction}}
Let $X$ be a $G$-space. 
\begin{itemize} 
\item Restriction: If $H \subset G$ then $X$ is naturally a $H$-space, and there is an induced map $H_G^*(pt) \to H_H^*(pt)$.
\[H_H^*(X)=H_H^*(pt) \otimes_{H_G^*(pt)} H_G^*(X)\]
\item Induction: If $G \subset K$ then $K\times_G X$ is naturally a $K$-space, and there is an induced map $H_K^*(pt) \to H_G^*(pt)$.  
\[H_K(K \times_G X)=H_G(X)\]
but as a $H_K(pt)$-module. 
\end{itemize}
\end{prop}

\begin{exa} Let $G=GL_n$. We have a left-right action of the upper Borel $B\subset GL_n$ on $G\times_B G$.  We compute $H_{B \times B}(G \times_B G)$ in the following steps: 
\[
H^*_{B\times B}(B)=H^*(B)=S,\]
so by induction
\[H^*_{G\times B}(G\times_B B)=S\in S^W \text{-mod-} S,\]
therefore by restriction
\[H^*_{B\times B}(G)=S \otimes_{S^W} S \in S \text{-mod-} S\]
and by induction again 
\[H^*_{G \times B}(G \times_B G)=S \otimes_{S^W} S \in S^W \text{-mod-} S\]
and by restriction
\[H^*_{B \times B}(G \times_B G)=S \otimes_{S^W} S \otimes_{S^W} S \in S \text{-mod-} S. \]
\end{exa}

\subsection{The equivariant DeRham model}

Let $G$ be a Lie group with Lie algebra $\mathfrak{g}$. For a smooth $G$-manifold $M$ we can define equivariant differential forms, for more details see \cite{getzlervergne}. The equivariant differential forms are differential form valued polynomial functions on $\mathfrak{g}$:
\[\Omega_G(M)=\left\{\alpha:\mathfrak{g}\to
\Omega(M):\alpha(gX)=g\alpha(X) \text{ for }g\in G , X \in \mathfrak{g}\right\} =(\C[\mathfrak{g}]\otimes \Omega(M))^{G}\]
where $(g\cdot \alpha)(X)=g\cdot (\alpha(g^{-1}\cdot X))$. Here $\C[\mathfrak{g}]$ denotes the algebra of complex values polynomial functions on $\mathfrak{g}$.
 
We define an equivariant exterior differential $d_G$ on  $\C[\mathfrak{g}]\otimes \Omega(M)$ by the formula 
\[(d_G\alpha)(X)=(d-\iota(X_M))\alpha(X),\]
where $\iota(X_M)$ denotes the contraction by the vector field $X_M$. 
This increases the degree of an equivariant form by one if the $\mathbb{Z}$-grading is given on $\C[\mathfrak{g}]\otimes \Omega(M)$ by
\[\deg(P \otimes \alpha)=2\deg(P)+\deg(\alpha)\]
for $P \in \C[\mathfrak{g}], \alpha \in \Omega(M)$. 
The homotopy formula $\iota(X)d+d\iota(X)=\mathcal{L}(X)$ implies that 
\[d_G^2(\alpha)(X)=-\mathcal{L}(X)\alpha(X)=0\]
for any $\alpha \in \C[\mathfrak{g}]\otimes \Omega(M)$, and therefore $(d_G,\Omega_G(M))$ is a complex.
\begin{defi}
The equivariant cohomology of the $G$-manifold $M$ is the cohomology of the complex $(d_G,\Omega_G(M))$:
\[H_G^*(M)=H^*_{d_G}\]
\end{defi}

Note that $\alpha \in \Omega_G(M)$ is equivariantly closed if
\[\alpha(X)=\alpha(X)_0+\ldots +\alpha(X)_n \text{ such that }
\iota(X_M)\alpha(X)_i=d\alpha(X)_{i-2}.\]

Here $\alpha(X)_i \in \Omega_i(M)$ is the degree-$i$ part of $\alpha(X) \in \Omega(M)$. In other words, $\alpha_i: \mathfrak{g} \to \Omega^i(m)$ is a polynomial function. 

The functoriality properties of equivariant cohomology now come for free:
\begin{enumerate}
\item If $H \to G$ is a homomorphism of Lie groups then the restriction map $\C[\mathfrak{g}] \to \C[\mathfrak{h}]$ induces a homomorphism of differential graded algebras $\Omega_G(M) \to \Omega_H(M)$ and finally a homomorphism $H_G(M) \to H_H(M)$. 
\item If $\phi:N \to M$ is a map of $G$-manifolds which interwines the actions of $G$ then pull-back by $\phi$ induces a homomorphism of differential graded algebras $\phi^*:\Omega_G(M) \to \Omega_G(N)$ and homomorphism $H_G(M) \to H_G(N)$. 
\end{enumerate} 

\section{Equivariant localisation}

\subsection{Integrating equivariant forms} If $G$ is a Lie group and $M$ is a $G$-manifold, we can integrate equivariant forms obtaining a map
\[\int_M : \Omega_G(M) \to \C[\mathfrak{g}]^G\]
 by the formula 
\[\left(\int_M \alpha\right)(X)=\int_M \alpha(X)=\int_M \alpha_{[n]}(X)\]

That is, if $\alpha$ is an equivariant differential form, then we integrate the top degree part of it, and a result is a polynomial function on $\mathfrak{g}$. This is well-defined: if $\alpha$ is equivariantly exact, i.e. $\alpha=d_G \beta$ for some $\beta \in \Omega_G(M)$ then $\alpha(X)_n=d \beta(X)_n$, and therefore $\int_M \alpha(X)=0$. Thus if $\alpha$ is equivariantly closed then $\int_M \alpha$ only depends on the equivariant cohomology class represented by $\alpha$. 

It can be shown (see Proposition 7.10 in \cite{getzlervergne}) that if $G$ is a compact Lie group, and $M_0(X)$ is the zero locus of the vector field $X_M$, then the form $\alpha(X)_n$ is exact outside $M_0(X)$. This suggests that the integral $\int_M \alpha(X)$ only depends on the restriction $\alpha(X)|_{M_0(X)}$.

Here we state the localisation thorem in the special case when $X_M$ has isolated zeros.  

\begin{thm}[Atiyah-Bott \cite{ab}, Berline-Vergne \cite{BV}] Let $G=T$ be a complex torus, $M$ a $T$-manifold, $\alpha \in \Omega_T(M)$. Then 
\[\int_M\alpha=(2\pi)^l\sum_{p\in
M^T}\frac{\alpha_0(p)}{\mathrm{Euler}^T(T_pM)}\] 

In other words:
\[\int_M\alpha(X)=(2\pi)^l\sum_{p\in
M^T}\frac{\alpha(X)_0(p)}{\prod_i \lambda_i}\] where $\lambda_i$
are the weights of the Lie action
\[X:\xi\in T_pM \to [X_M(p),\xi]\in T_pM.\]
\end{thm}

Most often we apply localisation to compute certain intersection numbers on $M$. 
My favorite example illustrating the strength of the localisation method is the following.

\subsection{How many lines intersect 2 given lines and go through a point in $\PP^3$?}\label{subsection:example}

We think points, lines and planes in $\PP^3$ as $1,2,3$-dimensional subspaces in $\C^4$. For $R \in \grass(3,\C^4), L \in \grass(1,\C^4)$ define
\[C_2(R)=\{V\in \grass(2,4): V\subset R\},\ C_1(L)=\{V \in \grass(2,4):L \subset V\}\]

Standard Schubert calculus says that $C_1(L)$ (resp $C_2(R)$) represents the cohomology class $c_1(\tau)$ (resp $c_2(\tau)$) where $\tau$ is the tautological rank 2 bundle over $\grass(2,4)$.

So the answer is  
\[C_1(L_1)\cap C_1(L_2) \cap C_2(R)=\int_{\grass(2,4)}c_1(\tau)^2c_2(\tau).\]
 
Apply equivariant localization. The sufficient data are the following.  
\begin{itemize}
\item The diagonal torus $T^4 \subset GL(4)$ acts on $\C^4$ with weights $\mu_1,\mu_2,\mu_3,\mu_4\in \mathfrak{t}^*\subset H_T^*(pt)$.  
\item The induced action on $\grass(2,4)$ has ${4 \choose 2}$ fixed points, the coordinate subspaces indexed by $(i,j)$. 
\item The tangent space of $\grass(2,4)$ at $(i,j)$ is $(\C^2)_{i,j}^* \otimes \C_{s,t}^2$, where $\{s,t\}=\{1,2,3,4\}\setminus \{i,j\}$, and $\C^2_{i,j}\in \grass(2,4)$ is the subspace spanned by the $i,j$ basis. Therefore, the weights on $T_{(i,j)}\grass$ are $\mu_s-\mu_i,\mu_s-\mu_j$ with $s\neq i,j$. 
\item The weights of $\tau$ are identified with the Chern roots, so $c_i(\tau)$ is represented by the $i$th elementary symmetric polynomial in the weights of $\tau$.  
\end{itemize}
ABBV localisation then gives
\begin{equation}\label{ABBV}
\int_{Gr(2,4)}c_1(\tau)^2c_2(\tau)=\sum_{\sigma \in S_4/S_2} \sigma \cdot \frac{(\mu_1+\mu_2)^2\mu_1\mu_2}{(\mu_3-\mu_1)(\mu_4-\mu_1)(\mu_3-\mu_2)(\mu_4-\mu_2)}=2.
\end{equation} 
On the right hand side we sum over all ${4 \choose 2}$ fixed points by taking appropriate permutation of the indices. 

It is not clear at first glance, why this rational expression is an integer. But it turns out that the sum is indepenent of $\mu_i$'s and it is 2. 

\subsection{Iterated residues}

We saw in the previous example that the ABBV localisation results a sum of rational expressions, but adding these together is not an obvious task. There is a short and elegant way to do this by identifying the summands as iterated residues of a certain meromorphic differential form on $\C^d$ for some $d$, and then by applying the Residue theorem saying that the sum of the residues at finite points is equal to minus the residue at infinity. 

The set-up is the following.

\begin{itemize}
\item $z_1,\ldots, z_d$ -- coordinates on $\C^d$.
\item $\omega_1,\dots,\omega_N$ -- affine linear forms on $\C^d$; $\omega_i=a_i^0+a_i^1z_1+\ldots + a_i^dz_d$. \item $h(\bz)$ a function $h(z_1\ldots z_d)$, and $\dbz=dz_1\wedge\dots\wedge dz_d$ holomorphic $d$-form.
\end{itemize}
\begin{defi} We define the iterated residue of $\frac{h(\bz)\,\dbz}{\prod_{i=1}^N\omega_i}$ at infinity as follows
\begin{equation}
  \label{defresinf}
 \res_{z_1=\infty}\ldots \res_{z_d=\infty}\frac{h(\bz)\,\dbz}{\prod_{i=1}^N\omega_i}
  \overset{\mathrm{def}}=\left(\frac 1{2\pi i}\right)^d
\int_{|z_1|=R_1}\ldots
\int_{|z_d|=R_d}\frac{h(\bz)\,\dbz}{\prod_{i=1}^N\omega_i},
 \end{equation}
 where $1\ll R_1 \ll \ldots \ll R_d$. The torus $\{|z_m|=R_m;\;m=1\ldots
 d\}$ is oriented in such a way that $\res_{z_1=\infty}\ldots
 \res_{z_d=\infty}\dbz/(z_1\cdots z_d)=(-1)^d$.
\end{defi}

In practice, the iterated residue \ref{defresinf} may be computed
using the following {\bf algorithm}: for each $i$, use the expansion
 \begin{equation}
   \label{omegaexp}
 \frac 1{\omega_i}=\sum_{j=0}^\infty(-1)^j\frac{(a^{0}_i+a^1_iz_1+\ldots
   +a_{i}^{q(i)-1}z_{q(i)-1})^j}{(a_i^{q(i)}z_{q(i)})^{j+1}},
   \end{equation}
   where $q(i)$ is the largest value of $m$ for which $a_i^m\neq 0$,
   then multiply the product of these expressions with $(-1)^dh(z_1\cdots
   z_d)$, and then take the coefficient of $z_1^{-1} \ldots z_d^{-1}$
   in the resulting Laurent series.

\begin{exa}
\begin{itemize}
\item $\frac{1}{z_1(z_1-z_2)}$ has two different Laurent expansions, but on $|z_1| \ll |z_2|$ we use $\frac{1}{z_1(z_1-z_2)}=\sum_{i=0}^\infty (-1)^i\frac{z_1^{i-1}}{z_2^{i+1}}$ to get $\res_\infty \frac{1}{z_1-z_2}=1$. 
\item $\res_{\bz=\infty}\frac{1}{(z_1-z_2)(2z_1-z_2)}=\mathrm{coeff}
_{(z_1z_2)^{-1}}\frac{1}{z_2^2}(1+\frac{z_1}{z_2}+\frac{z_1^2}{z_2^2}+\ldots)
(1+\frac{2z_1}{z_2}+\frac{4z_1^2}{z_2^2}+\ldots)=3$
\end{itemize}

\end{exa}

Let's turn back to our toy example presented in \S \ref{subsection:example}. Define the differential form 
\[\omega=\frac{(z_2-z_1)^2(z_1+z_2)^2z_1z_2\dbz}{\prod_{i=1}^4(\mu_i-z_1)\prod_{i=1}^4(\mu_i-z_2)}\]
This is a meromorphic form in $z_2$ on $\PP^1$ with poles at $z_2=\mu_i, 1\le i\le 4$ and $z_2=\infty$. The poles at $\mu_i$ are non-degenerate and therefore applying the Residue Theorem we get 
\[\res_{z_2=\infty}\omega=\sum_{i=1}^4\underbrace
{-\frac{(\mu_i-z_1)^2(\mu_i+z_1)^2\mu_iz_1dz_1}{
\prod_{j=1}^4(\mu_j-z_1)\prod_{j\neq i}(\mu_j-\mu_i)}}_{z_2=\mu_i}=\sum_{i=1}^4
-\frac{(\mu_i-z_1)(\mu_i+z_1)^2\mu_iz_1dz_1}{
\prod_{j\neq i}(\mu_j-z_1)\prod_{j\neq i}(\mu_j-\mu_i)}\]
Doing the same again with $z_1$ we get
\[\res_{z_1=\infty}\res_{z_2=\infty}\omega=\sum_{i=1}^4 \sum_{j\neq i}
-\frac{(\mu_i-\mu_j)(\mu_i+\mu_j)^2\mu_i\mu_j}{
\prod_{k\neq i,j}(\mu_k-\mu_j)\prod_{j\neq i}(\mu_j-\mu_i)}=\]
\[=\sum_{i=1}^4 \sum_{j\neq i}
\frac{(\mu_i+\mu_j)^2\mu_i\mu_j}{
\prod_{k\neq i,j}(\mu_k-\mu_j)\prod_{k\neq i,j}(\mu_k-\mu_i)}=\int_{Gr(2,4)}c_1(\tau)^2c_2(\tau).\]
On the other hand, using the above algorithm by expanding the rational form $\omega$ we get
\[\res_{z_1=\infty}\res_{z_2=\infty}\omega=2.\]

We give an other example, the so called Giembelli-Thom-Porteous formula in section \ref{giambelli}.

\subsection{Localisation on partial flag manifolds}

The following Proposition is a far-reaching generalization of the idea presented in the previous section, and it provides a meromorphic differential form whose residue at infinity gives back the localisation formula for a large class of forms. 

Let \[\flag_d(n)=\{V_1\subset \ldots \subset V_d\subset \C^n:\dim(V_i)=i\}\]
denote the full flags of $d$-dimensional subspaces of $\C^n$. The maximal torus $T \subset GL(n)$ acts on $\flag_d(n)$, and the fixed points are parametrized by coordinate flags corresponding to certain permutations $\sigma \in (\flag_d(n))^T$. The Chern classes of the tautological rank-$d$ bundle over $\flag_d(n)$ are elementary symmetric polynomials in the weight of $T$ on $\C^n$, and the intersection numbers of $\tau$ can be computed as iterated residues according to  
\begin{prop}[\cite{bsz}]\label{propflag}
  \label{flagresidue} Let $Q(\bz)=Q(z_1,\ldots, z_d)$ be a polynomial on $\C^d$ of degree $\dim(\flag_d(n))$.  Then 
\begin{equation}
  \label{flagres}
\sum_{\sigma \in (\flag_d(n))^T}
\frac{Q(\lambda_{\sigma\cdot 1}\cdots\lambda_{\sigma\cdot d})}
{\prod_{1\leq m\leq d}\prod_{i=m+1}^n(\lambda_{\sigma\cdot
    i}-\lambda_{\sigma\cdot m})}=\sires
\frac{\prod_{1\leq m<l\leq d}(z_m-z_l)\,Q(\bz)\dbz}
{\prod_{l=1}^d\prod_{i=1}^n(\lambda_i-z_l)}
\end{equation}
where the permutation $\sigma=(\sigma(1),\ldots ,\sigma(n))\in \flag_d(n))^T=S_n/S_{n-d}$ represents the torus fixed flag $\C e_{\sigma(1)}\subset \ldots \subset \C e_{\sigma(1)} \oplus \ldots \oplus \C e_{\sigma(d)} \subset \C^n$.
\end{prop}

\section{Singularities of maps}

The first problem we address goes back to the 1950's and the work of Ren\'e Thom \cite{thom}. For more details of the history and background of the problem see \cite{arnold,bsz,damon, FR2, FR3, kazarian}.

The usual set-up for studying singularities of map germs is the following.

\textbf{Set up:} We fix integers $k\le n\le m$.
\begin{itemize}
 \item Let $A$ be a nilpotent algebra, $\dim A/\C=k$. We will take $A_k=z\C[z]/z^{k+1}$.
 \item Define $J_k(n,m)=\left\{p=(p_1,\ldots, p_m)\in
\mathrm{Poly}(\C^n,\C^m):\deg p_i \le k,p_i(0)=0 \right\}$. This is the vector space of $k$-jets of map germs $f:(\C^n,0) \to (\C^m,0)$.
 \item Let $\Sigma_A=\left\{p \in J_k(n,m):\C[x_1,\ldots, x_n]/\langle p_1,\ldots, p_m \rangle
=A \right\}$ be the set of map-germs with local algebra isomorphic to $A$.
 \item The germs $J_k^\reg(n,n)$ with non-degenerate linear part form a group, and $J_k^\reg(n,n) \times J_k^\reg(m,m)$ naturally acts on $J_k(n,m)$ with
\[(A,B)p=BpA^{-1}\]
These are the polynomial reparametrisations of map germs.
\end{itemize}


A central problem in global singularity theory is the computation of the cohohomology classes of singularity loci of holomorphic maps between complex manifolds. Given a holomorphic map $f:N^n \to M^m$
define 
\[Z(f)=\left\{p \in N| \widehat{f}_p\in
\Sigma_A\right\},\]
where $\widehat{f}_p$ is the germ of $f$ at $p\in N$.

It was already known by Thom for real manifolds and differentiable maps between them, which later was extended to the complex case in \cite{damon}---and is now known as the Thom principle---that for generic map $f$, $Z(f)$ represents a cohomology cycle and there is a well-defined polynomial
\[MD^{n \to m}_A \in \C[x_1,\ldots ,x_n,y_1,\ldots ,y_m]^{S_n\times S_m}\]
such that  
\[[Z(f)]=MD_A(TN,f^*(TM))\in H^*(N,\C).\]
Here $MD_A$ stands for multidegree, for explanation see the next section. 

Furthermore, in \cite{HK} Haefliger and Kosinski proves that if 
\[c(q)=c_0+c_1q+c_2q^2+\ldots =\frac{c(f^*(TM))}{c(TN)}=
\frac{\prod_{m=1}^k(1+\theta_mq)}{\prod_{i=1}^n(1+\lambda_iq)}\]
is the Chern classes of the difference bundle then  
\[MD_A(TN,f^*(TM))=\mathrm{Tp}_A^{k\to n}(c_1,c_2,\ldots)\]
That is, $MD_A$ is a polynomial in these difference Chern classes, and $\mathrm{Tp}_A$ is called the Thom polynomial of the algebra $A$. 

\subsection{Multidegrees}
The polynomial $MD_A$ stands for multidegree, which is also called equivariant Hilbert polynomial or equivariant Poincar\'e dual in the literature. This is defined for any $G$-invariant
subvarieties of a complex $G$-vector spaces (i.e. $G$-representations, where $G$ is a Lie algebra) as follows.

\noindent\textbf{Set up:}
\begin{enumerate}  \item $V=\C^N$ complex vector space, with
a $G$-action.
 \item $\Sigma \subset V$ is a $G$-invariant closed subvariety.
 \item $H_G^*(V)=H_G^*(pt)$ is the $G$-equivariant cohomology ring of
$V$. Recall that $H_{GL(d)}^*(pt)=\C[x_1,\ldots, x_d]^{S_d}$.
\end{enumerate}



We give two definitions of a polynomial 
\[\mdeg{\Sigma,V}\in H_G^{\mathrm{codim}(\Sigma \in
V)}(pt),\]
called the multidegree of $\Sigma$: one topological and one algebraic definition.

\textbf{Vergne's integral definition-topology}\\
If $\Sigma \subset V$ is a subvariety then $EG \times_G \Sigma\subset EG \times_G V$ represents a homology
cycle, and the multidegree is the ordinary Poincar\'e dual of the Borel construction $EG \times_G \Sigma$:
\[\mdeg{\Sigma,V}=PD(EG \times_G \Sigma\subset EG \times_G V).\]
By definition $\mdeg{\Sigma,V}\in H^*(EG \times_G V)=H_G^*(pt)$ is a polynomial. 

\begin{thm}[\cite{vergne}] There is an equivariant Thom class:
\[\thom_G(V)\in H_G^{\dim V}(V)\]
such that for any $\Sigma \subset V$ $G$-invariant subvariety
\[\mdeg{\Sigma,V}=\int_{\Sigma}\thom_G(V).\]
\end{thm}



We give an other, more algebraic definition of the multidegree, which also provides an algorith to compute these polynomials.\\

\noindent \textbf{Axiomatic definition}

\begin{thm}[\cite{milsturm}] Let $\Sigma \subset V$ be a $G$-invariant subset of the $G$-representation $V$. Then $\mdeg{\Sigma,V}$ is characterized by the following axioms:
\begin{description}
 \item[additivity] If $\Sigma \doteq \cup \Sigma_i$ is the set of maximal irreducible components of $\Sigma$, then 
\[\mdeg{\Sigma,V}=\sum_{i=1}^c \mathrm{mult}(\Sigma_i)
\cdot \mdeg{\Sigma_i,W}.\]
 \item[degeneration] The multidegree is constant under flat deformation of $\Sigma$.
 \item[normalization] For $T$-invariant linear subspaces of $V$, $\mdeg{\Sigma,V}$ is defined to be equal to the product of weights in the normal
  direction.
\end{description}
\end{thm}

The recipe to compute the multidegree (although this recipe often ends up with difficult commutative algebra computations) is to choose a proper flat deformation of $\Sigma$ into the union of coordinate spaces, that is, to deform its ideal into a monomial ideal. For example, choosing a monomial order on the coordinate ring, the initial ideal is monomial, and then normalization and additive properties of the multidegree give you the result.  

\begin{exa} $(\C^*)^3$ acts on $\C^4$ with weights
$\eta_1,\ldots ,\eta_4$. Let $\eta_1+\eta_2=\eta_3+\eta_4$, and
\[\Sigma=\mathrm{Spec}(\C[y_1,y_2,y_3,y_4]/(y_1y_2-y_3y_4)).\]
Define the flat deformation 
\[\Sigma_t=\mathrm{Spec}(\C[y_1,y_2,y_3,y_4]/(y_1y_2-ty_3y_4)),\]
For $t=0$, $\Sigma_0=\left\{y_1y_2=0 \right\}$, so normalization and additivity says
\[\mdeg{\Sigma,\C^4}=\eta_1+\eta_2=\eta_3+\eta_4.\]
\end{exa}

Now we can state Thom's principle more precisely:
\begin{thm}[\cite{damon,thom}]\label{Thom} Let $\Sigma_A \subset J_k(n,m)$ denote the set of germs with local algebra isomorphic to $A$. This is a $\mathrm{GL}_n \times \mathrm{GL}_m$-invariant subvariety of $J_k(n,m)$, and 
\[MD_A^{n\to m}=\mathrm{mdeg}^{\mathrm{GL}_n \times \mathrm{GL}_m}[\Sigma_A,J_k(n,m)].\]
\end{thm}

\subsection{Degeneraci loci of sections via localization}\label{giambelli}

Here is an other illustrating example for transformation of localisation formulas into iterated residues. Given a a rank-$n$ vector bundle on a manifold $M$, and $n$ generic sections $\sigma_1,\ldots, \sigma_n$, it is an old question in topology to determine the cohomology class dual to the locus where the sections are linearly dependent. This class is the Thom polynomial $\mathrm{Tp}_A$ with $A=t\C[t]/t^2$, and we have
\[\Sigma_A=\Sigma_1 = \{A\in \Hom(n,m);\; \dim \ker A= 1\}=
\{A\in \Hom(n,m)\exists! [v]\in \PP^{n-1}:\,Av=0\}.\]
The goal is to compute $\mdeg{\Sigma_1,\Hom(n,m)}$. 

We have the fibration $\pi:\Sigma_1 \rightarrow \PP^{n-1}$ sending a linear map to its kernel. This is equivariant with respect to the $GL_n \times GL_m$ action, and $GL_m$ acts fibrewise whereas $GL_n$ acts on $\PP^{n-1}$. According to Vergne's definition, we want to integrate the equivariant Thom class over $\Sigma_1$. The idea is to integrate first over the base $\PP^{n-1}$ and then along the fibres, and to apply ABBV localisation on $\PP^{n-1}$. 

We have $n$ fixed points on $\PP^{n-1}$ corresponding to the coordinate axes. Let $\lambda_1,\ldots, \lambda_n$ denote the weights of $T^n \subset GL_n$ on $\C^n$. The weights of
$T_{p_i}\PP^{n-1}$ are $\{\lambda_s-\lambda_i;\;s\neq i\}$, and the 
fibre at $p_i$ is the set of matrices $A$ with all entries in the
$i$th column vanishing. The normalization axiom says that the multidegree of the fibre at $p_i$ is
$\prod_{j=1}^m(\theta_j-\lambda_i)$, so:

\[\mdeg{\Sigma_1,\Hom(
n,k)}=\int_{\Sigma_1}\mathrm{Thom}_{(\C^*)^{n+m}}(\hm
nm)=\]
\[=\int_{\PP^{n-1}}\int_{fibre}\mathrm{Thom}_{(\C^*)^{n+m}}=
\sum_{i=1}^n\frac{\prod_{j=1}^m(\theta_j-\lambda_i)}{\prod_{s\neq
i}(\lambda_s-\lambda_i)}
\]

Consider the rational differential form
\[
-\frac{\prod_{j=1}^{m}(\theta_j-z)}
{\prod_{i=1}^{n}(\lambda_i-z)}\,dz.\] The residues of this form at
finite poles: $\{z=\lambda_i;\;i=1\ldots n\}$ exactly recover the
terms of the above sum. Applying the residue theorem, and change of variables
$z=-1/q$, we get

\[\mdeg{\Sigma_1,\hm nm}=\mathrm{res}_{q=0}\frac{\prod_{j=1}^{m}(1+q\theta_j)}
{\prod_{i=1}^{n}(1+q\lambda_i)}\,\frac{dq}{q^{m-n+2}}=c_{m-n+1},\]
where $c_{m-n+1}$ is the $m-n+1$th Chern class of the difference bundle $f^*(TM)-TN$. This gives us the Thom polynomial $Tp_{t\C[t]/t^2}^{n \to m}=c_{m-n+1}$. Note that it depends only on $m-n$.

\section{Computing multidegrees of singularities}

Recall the following notations from the previous section. 
\begin{itemize}
\item $J_k(n,m)=\left\{(p_1,\ldots, p_m)\in
\mathrm{Poly}(\C^n,\C^m):\deg p_i \le k,p_i(0)=0 \right\}$ is the set of $k$-jets of map germs.
\item $\Sigma_k=\left\{p \in J_k(n,m):\C[x_1,\ldots, x_n]/\langle p_1,\ldots, p_m \rangle
\cong z\C[z]/z^{k+1} \right\}$ the set of germs with $A_k$-singularity. 
\item $\mathcal{D}=J_k^\reg(n,n) \times J_k^\reg(m,m)$ naturally acts on $J_k(n,m)$ with $(A,B)p=BpA^{-1}$. Note that $\mathrm{GL}_n \times \mathrm{GL}_m \subset \mathcal{D}$
\end{itemize}

The goal now is to compute 
\[\tp_k^{n \to m}=\mathrm{mdeg}^{\mathrm{GL}_n \times \mathrm{GL}_m}[\Sigma_k,J_k(n,m)],\]
the Thom polynomial of Morin singularities. 


\begin{thm}[\textbf{The test curve model of Morin singularities},\cite{gaffney}]
\[\Sigma_k(n,m)\doteq \left\{\Psi \in J_k(n,m)|\exists \gamma
\in J_k^\reg(1,n) \text{ such that\ } \Psi \circ \gamma =0 \text{ in } J_k(1,m) \right\}.\]
Here $\doteq$ denotes birational equality, that is their Zariski closure are equal.   
\end{thm}
\begin{diagram}[LaTeXeqno,labelstyle=\textstyle]\nonumber
 (\C,0) & \rTo^\gamma & (\C^n,0) & \rTo^{\Psi} & (\C^m,0)
\end{diagram}
Note that if $\varphi \in J_k^{\reg}(1,1)=\mathbf{G}_k$,
then
\[\Psi \circ \gamma=0\ \ \Rightarrow \ \ \Psi \circ (\gamma \circ \varphi)=0 \]
\begin{diagram}[LaTeXeqno,labelstyle=\textstyle]
  (\C,0) & \rTo^\varphi & (\C,0) & \rTo^\gamma & (\C^n,0) & \rTo^{\Psi} & (\C^m,0)
\end{diagram}

It can be shown that for $\Psi \in J_k(n,m)$ whose linear part has corank 1 
\[\Psi \circ \gamma_1=\Psi \circ \gamma_2=0 \Leftrightarrow \exists \alpha \in J_k^\reg(1,1) \text{ s.t }\gamma_1 =\gamma_2 \circ \alpha.\]
Therefore: 

\begin{prop}\label{fibres} The Zariski open subset $\Sigma^0_k=\left\{\Psi \in \Sigma_k:\dim \ker \Psi=1
\right\}\doteq \Sigma_k$ fibres with linear fibres over
$J_k^{reg}(1,n)/\mathbf{G}_k$. 
\end{prop}

What are these fibres, and why are they linear? If $\gamma=v_1t+v_2t^2+\ldots +v_dt^d \in J_k^\reg(1,n)$ with $v_i \in \C^n$ and $v_1 \neq 0$ and 
$\Psi(v)=Av+Bv^2+\ldots$ with $A\in \mathrm{Hom}(\C^n,C^k)$, $B\in
\mathrm{Hom}(\mathrm{Sym}^2(\C^n),\C^k),etc$, 
then $\Psi \circ \gamma=0$ is equivalent with the following $k$ equations:
\begin{eqnarray}
& A(v_1)=0, \\ \nonumber & A(v_2)+B(v_1,v_1)=0, \\ \nonumber
& A(v_3)+2B(v_1,v_2)+C(v_1,v_1,v_1)=0, \\
& ... \nonumber
\end{eqnarray}

For fixed $\gamma=(\gamma_1,\ldots, \gamma_k)$ these are linear equations determining the fibre. 
According to Proposition \ref{fibres}
\[\Sigma_k(n,m)\doteq \bigcup\left\{\mathrm{Sol}_\gamma \mid \gamma \in J_k^{reg}(1,n) \right\},\]
where 
\[\mathrm{Sol}_\gamma=\ann(\gamma)\otimes \C^m \subset J_k(n,m)\]
is the annihilator tensored by $\C^k$.

To linearize the action of $\mathbf{G}_k$ let's make the following identifications
\begin{itemize}
\item Identify  $J_k(1,n)$ with $\Hom(\C^k,\C^n)$ by putting the coordinates $\gamma=(v_1,\ldots, v_k)$ into the columns of a matrix;
\item Identify $J_k(n,1)$ with $\mathrm{Sym}^{\le k}\C^n=\oplus_{i=1}^k\mathrm{Sym}^i\C^n$, and then $J_k(n,m)=\mathrm{Sym}^{\le k}\C^n \otimes \C^m$.   
\end{itemize}

Then $\mathbf{G}_k$ acts on $J_k(1,n)$ by multiplication on the right by the following matrix group:
\begin{equation}\label{group}\left\{\left(
\begin{array}{ccccc}
\alpha_1 & \alpha_2   & \alpha_3          & \ldots & \alpha_k\\
0        & \alpha_1^2 & 2\alpha_1\alpha_2 & \ldots & 2\alpha_1\alpha_{k-1}+\ldots \\
0        & 0          & \alpha_1^3        & \ldots & 3\alpha_1^2\alpha_ {k-2}+ \ldots \\
0        & 0          & 0                 & \ldots & \cdot \\
\cdot    & \cdot   & \cdot    & \ldots & \alpha_1^k
\end{array}
 \right): \alpha_1 \in \C^* , \alpha_i \in \C  \right\};
\end{equation}
where the polynomial in the $(i,j)$ entry is
\[p_{i,j}(\bar{\alpha})=\sum_{a_1+a_2+\ldots +a_i=j}\alpha_{a_1}\alpha_{a_2} \ldots \alpha_{a_i}.\]

This group is the central object of our study in this paper. It is a non-reductive linear subgroup of $GL_k$, and therefore Mumford's geometric invariant theory \cite{git} does not help us in handling the quotient $J_k^\reg(1,n)/\mathbf{G}_k$. The following construction, which was the starting point of a general construction in \cite{bk2} first appeared in \cite{bsz}.

Define the map
\begin{equation}\label{rho} \rho:\mathrm{Hom}(\C^k,\C^n) \rightarrow
\mathrm{Hom}(\C^k,\mathrm{Sym}^{\le k}\C^n)
\end{equation}
\[\rho(v_1,\ldots, v_k)=(v_1,v_2+v_1^2,\ldots, \sum_{a_1+a_2+\ldots +a_i=k}v_{a_1}v_{a_2} \ldots v_{a_i}),\]
where in the $j$th coordinate we sum over all ordered partitions of $j$ into positive integers. Note that these correspond to the monomials in $j$th column of the matrix $\mathbf{G}_k$. For more details see \cite{bk2}.
\begin{thm}[\textbf{\cite{bsz}}]\label{bercziszenes}
Let $\mathrm{Hom}^0(\C^k,\C^n)=\{(v_1,\ldots, v_k \in \Hom(\C^k,\C^n):v_1 \neq 0\}=J_k^{reg}(1,n)$. Then $\rho$ (defined in \eqref{rho}) descends to an injective map on the orbits
\[\rho^{\grass}: \mathrm{Hom}^{0}(\C^k,\C^n)/\mathbf{G}_k \hookrightarrow
\grass(k,\mathrm{Sym}^{\le k}\C^n).\]
and therefore descends also to 
\[\rho^{\flag}: \mathrm{Hom}^0(\C^k,\C^n)/\mathbf{G}_k \hookrightarrow
\flag_k(\mathrm{Sym}^{\le k}\C^n)\]
\end{thm}
Composing with the Plucker embedding we get
\[\rho^{\proj}=Pluck \circ \rho^{\grass}: \mathrm{Hom}^0(\C^k,\C^n)/\mathbf{G}_k \hookrightarrow
\PP(\wedge^k(\mathrm{Sym}^{\le k}\C^n))\]


Note that $\rho$ is $GL_n$-equivariant with respect to the multiplication on the left on $\mathrm{Hom}^0(\C^k,\C^n)/\mathbf{G}_k$ and the induced action on $\grass(n,\symk)$ coming from the standard action on $\C^n$.  



This embedding allows us to give a geometric description of some generators in the invariant ring \[\C[J_k^{reg}(1,n)]^{U_k} \subset \C[f',\ldots, f^{(k)}],\]
where $U_k \subset \mathbf{G}_k$ is the maximal unipotent subgroup. Namely, the coordinate ring of the image is a subring of the invariant ring:
\[\C[\mathrm{im}(\rho)]\subset \C[J_k^{reg}(1,n)]^{U_k}\]

$\C[J_k^{reg}(1,n)]^{U_k}$ has been long studied in relation with the Green-Griffiths conjecture. In his seminal paper \cite{dem}, Demailly suggested a strategy to the Green-Griffiths conjecture through the investigation of the invariant jet differentials. These are sections of a bundle, whose fibres are canonically isomorphic to the invariant ring $\C[J_k^{reg}(1,n)]^{U_k}$, and he conjectured that   
this invariant ring is finitely generated, which was finally proved in \cite{bk} and \cite{bk2}, where we give give geometric description of a generating set. Note that $\mathbf{G}_k$ is a non-reductive group, and therefore classical Geometric Invariant Theory (\cite{git}) to describe $\mathrm{Spec}(\C[J_k^{reg}(1,n)]^{U_k})$ does not apply. For an introduction on non-reductive group actions see \cite{dk}.  

Following Demailly's notation, let $(f',f'',\ldots,f^{(k)}) \in J_k^\reg(1,n)$ denote the $k$-jet of a germ $f$ and $f_i^{(j)}$ denote the $i$th coordinate of the $j$th derivative, $1\le i \le n, 1\le j \le k$. This is a simple rescaling, namely $v_i=f^{(i)}/i!$.

\begin{thm}[\cite{bk}] 
Let 
\[I=(\Delta_{\bi_1,\ldots, \bi_k}(f):(\bi_1,\ldots, \bi_k) \in {\dim(\mathrm{Sym}^{\le k}(n)) \choose k}) \vartriangleleft \C[f',\ldots,f^{(k)}]\]
be the ideal generated by the $k \times k$ minors of 
$\rho(f'\ldots, f^{(k)})\in \Hom(\C^k,\mathrm{Sym}^{\le k}\C^n )$. 
Then 
\[I \subset \C[J_k^{reg}(1,n)]^{U_k}\]
\end{thm}

\begin{exa} $n=2,k=4$. In this case
\[J_4^{\mathrm{reg}}(1,2)=\{(f_1',f_2',f_1'',f_2'',f_1''',
f_2''',f_1'''',f_2'''')\in
(\C^2)^4;(f_1',f_2')\neq (0,0)\},\] and fixing a basis $\{ e_1,e_2 \}$ of
$\C^2$ and
\[ \{e_1,e_2,e_1^2,e_1e_2,e_2^2,e_1^3, \ldots, e_1e_2^4, e_2^4\} \]
of $\Sym^{\le 4}\C^2$ the map $\rho:J_4(1,2) \to \Hom(\C^4,\Sym^{\leq 4}\C^2)$ sends 
\[ (f_1',f_2',f_1'',f_2'',f_1''',f_2''',f_1'''',f_2'''')\] 
to a $4 \times 15$ matrix, whose first $5$ columns (corresponding to $\Sym^{\le 2} \C^2$) are 
\[
\left(
\begin{array}{ccccc}
f_1' & f_2' & 0 & 0 & 0   \\
\frac{1}{2!}f_1'' & \frac{1}{2!}f_2'' & (f_1')^2 & f_1'f_2' & (f_2')^2  \\
\frac{1}{3!}f_1''' & \frac{1}{3!}f_2''' &  f_1'f_1'' &
(f_1'f_2''+f_1''f_2') & f_2'f_2''   \\
\frac{1}{4!}f_1'''' & \frac{1}{4!}f_2'''' & \frac{2}{3!}f_1'f_1'''+\frac{1}{2!2!}(f_1'')^2 & \frac{2}{3!}(f_1'f_2'''+f_1'''f_2')+\frac{1}{2!}f_1''f_2''&  \frac{2}{3!}f_2'f_2'''+\frac{1}{2!2!}(f_2'')^2  \\
 \end{array}
\right),\]
and next four columns (corresponding to $\Sym^3 \C^2$) are
\[\left(\begin{array}{cccc}
  0 & 0 & 0 & 0  \\
  0 & 0 & 0 & 0  \\
 (f_1')^3 & (f_1')^2f_2' & f_1'(f_2')^2 & (f_2')^3 \\
 \frac{3}{2!}((f_1')^2f_1'') & \frac{3}{2!}((f_1')^2f_2''+2f_1'f_2'f_1'')& \frac{3}{2!}((f_2')^2f_1''+2f_2'f_1'f_2'')& \frac{3}{2!}((f_2')^2f_2'')  
\end{array}\right)
,\]
and the remaining five columns (corresponding to $\Sym^3 \C^3$) are
\[\left(\begin{array}{ccccc}
0 & 0 & 0 & 0 & 0 \\
0 & 0 & 0 & 0 & 0 \\
0 & 0 & 0 & 0 & 0 \\
(f_1')^4 & (f_1')^3f_2' & (f_1')^2(f_2')^2 & f_1'(f_2')^3 & (f_2')^4
\end{array} \right).\]

Then the weight $1+2+3+4=10$ piece $\C[J_4(1,2)]^{U_4}_{10}$ of the invariant algebra 
$\C[J_4(1,2)]^{U_4}$
is generated by the $4\times 4$ minors of this $4 \times 15$ matrix.

\end{exa}

\subsection{The computation: double localization+vanishing theorem}

According to Proposition \ref{fibres} and Theorem \ref{bercziszenes} we have the following picture, also called the snowman-model after the figure in Section \S 6 in \cite{bsz}:
\begin{diagram}[LaTeXeqno,labelstyle=\textstyle]\label{model}
\Sigma_k(n,m)  & \subset J_k(n,m) \\
 \dTo & \\
 J_k^{\reg}(1,n)/\mathbf{G}_k & \subset \mathrm{Flag}_k(\mathrm{Sym}^{\le k}\C^n) \\
\dTo & \\
 J_k^{\reg}(1,n)/B_k=\mathrm{Flag}_k(\C^n) &  
\end{diagram}
Here $B_k \subset GL_k$ is the upper Borel subgroup. 
Now we apply ABBV localisation to compute $\mdeg{\Sigma_k(n,m),J_k(n,m)}$. According to Vergne \cite{vergne}, we have to compute 
\[\int_{\Sigma_k(n,m)}\thom(J_k(n,m)),\] 
and we do this in two steps: first we localise on $\flag_k(\C^n)$ and use Proposition \ref{propflag} to turn the localisation formula into an iterated residue. Then we integrate along the fibres. The fibres are canonically isomorphic to $B_k/\mathbf{G}_k$ and in the second step we apply ABBV localisation on the image $\overline{\rho(fibre)}\subset \grass(k,\symk)$. Surprisingly---for some unclear geometric reason---in this second localisation all fixed points but a distinguished one contributes $0$ to the sum, and a lenghty computation leads us in \cite{bsz} to


\begin{thm}[\cite{bsz}]
\begin{equation}\label{main}
\tp_k^{m-n}= \mathrm{Res}_{\mathbf{z}=\infty}
\frac{\prod_{i<j}(z_i-z_j)\,\mathcal{Q}_k(z_1\ldots z_k)}
{\prod_{i+j\le l \le k}(z_i+z_j-z_l)}\cdot
\prod_{l=1}^k c\left(\frac1{z_l}\right)\,z_l^{m-n}\;dz_l,
\end{equation}
where
\begin{itemize}
\item
We integrate on the cycle $|z_1|>|z_2|>\ldots |z_k|$, which determines the Laurent expansion.
 \item $c(q)=1+c_1q+c_2q^2+\ldots$
\item $Q_k(z_1,\ldots, z_k)$ is the multidegree of a Borel-orbit in $(\C^k)^* \otimes \Sym^2(\C^k)$, for details see \cite{bsz}, and   
\[Q_1=Q_2=Q_3=1,Q_4=2z_1+z_2-z_4\]
\end{itemize}
\end{thm}

The polynomial $Q_k$ is known up to $k \le 6$, but with enough computer capacity---in principle---it can be computed for any $k$. But no general formula is known at the moment. 

We give a concise--and not complete--summary of the history of Thom polynomial computations.
For a more detailed overview see \cite{bsz,kazarian}. 

\begin{itemize}
\item Multidegrees of singularities have been studied for nearly 60 years now. We call these Thom polynomials after the pioneering work of Ren\' e Thom in the 1950's. He proved the existence of these polynomials (\cite{thom}) for real manifolds and singularities of differentiable maps between them. Later Damon in \cite{damon} studied complex contact singularities.  
\item The case $k=1$ is the
classical formula of Porteous: $\tp_1^{n\to m}=c_{m-n+1}$. The $k=2$ case was computed by Ronga in \cite{ronga}. An explicit
  formula for $\tp_3^{n \to k}$ was proposed in \cite{BFR} and P. Pragacz has given a proof
  \cite{prag}. He also studied Thom polynomials in \cite{prag2,lp}.  Finally, using his method of restriction equations,
  Rim\'anyi \cite{rimanyi} was able to treat the $n=k$ case, and
  computed $\tp_k^{n\to n}$ for $k\le 8$ (cf. \cite{gaffney} for the
  case $d=4$). 
\item More recently, Kazarian (\cite{kazarian2}) has worked out a framework for computing Thom polynomials of contact singularities in general. He suggests studying certain non-commutative associative algebras to get a polynomial $Q_A$ and an iterated residue formula similar to \eqref{main} for any local algebra $A$. Unfortunately, the explicit computation of $Q_A$ is difficult, his description does not give more information for Morin singularities, where $Q_k$ is unknown for $k>6$. Rim\'anyi and Feh\'er in \cite{rf} compute Thom series for further singularity classes.  
The structure of Thom polynomials of contact singularities was also studied in \cite{FR2,FR3}.
\end{itemize}

Finally, let us mention a conjecture of R. Rim\'anyi about the positivity of these Thom polynomials. 

\begin{conjecture}[\textbf{Rim\'anyi} \cite{rimanyi}]\label{rimanyiconj}
\[\mathrm{Tp}_k^{m-n} \in \mathbb{N}[c_1,\ldots, c_{k(m-n+1)}],\]
i.e. the coefficients of the Thom polynomials are nonnegative. This would follow from the more general conjecture, that   \[\frac{\prod_{i<j}(z_i-z_j)\,\mathcal{Q}_k(z_1\ldots z_k)}
{\prod_{i+j\le l \le k}(z_i+z_j-z_l)}>0,\]
the coefficients of the Thom series are nonnegative. 
\end{conjecture}

\section{The Green-Griffiths conjecture}

First we list some results related to hyperbolic varieties and the Green-Griffiths conjecture.  
This is a selection of classical results and it is far from being complete. 

\subsection{Hyperbolic varieties}

Let $X$ be a complex manifold, $n=\dim_\C(X)$. $X$ is said to be hyperbolic 
\begin{itemize}
\item in the sense of Brody, if there
are no non-constant entire holomorphic curves $f: \C \to X$. 
\item in the sense of Kobayashi, if the Kobayashi-Royden pseudo-metric on $T_X$ is nondegenerate. This pseudo-metric is defined as follows. The infinitesimal Kobayashi-Royden metric is
\[k_X(\xi)=\inf\{\lambda>0:\exists f:\Delta \to X, f(0)=x,\lambda f'(0)=\xi\} \text{ for } \ x\in X, \xi \in T_{X,x}.\]
The Kobayashi pseudo-distance $d(x,y)$ is the geodesic pseudo-distance
obtained by integrating the Kobayashi-Royden infinitesimal metric. $X$ is hyperbolic in the sense of Kobayashi if $d(x,y)>0$ for $x\neq y$. 
\end{itemize}

The following theorem of Kobayashi tells that positivity of the cotangent bundle implies hyperbolicity.


\begin{thm}[\textbf{Kobayashi, '75}] $X$-smooth projective variety with ample cotangent bundle. Then $X$ is hyperbolic.
\end{thm}

Conversely, 

\begin{conjecture}
If a compact manifold $X$ is
hyperbolic, then it should be of general type, i.e. $K_X=\wedge^nT^*X$
should be big. (That is, $X$ has maximal Kodaira dimension, i.e. $\dim \oplus_{i=0}^\infty H^0(X,K^i)=\dim X$.)
\end{conjecture}

\begin{conjecture}[\textbf{Green-Griffiths, '79}] Let X be a
projective variety of general type. Then there exists an
algebraic variety $Y \varsubsetneq X$ such that for all non-constant
holomorphic $f : \C \to X$ one has $f(\C)\subset Y$.
\end{conjecture}

\textbf{Diophantine properties}
\begin{thm}[\textbf{Faltings, '83}] A curve of genus greater than 1 has only finitely many rational points. 
\end{thm}

\begin{thm}[\textbf{Moriwaki, '95}] Let $K$ be a number field (finitely generated over $\Q$), and $X$ a smooth projective variety. If $T^*X$ is ample and globally generated then $X(K)$ is finite. 
\end{thm}

\begin{conjecture}[\textbf{Lang, '86}] \begin{enumerate}
\item If a projective variety $X$ is hyperbolic, then it is mordellic, i.e. $X(K)$ is finite for any $K$ finitely generated over $\Q$.
\item Let $Exc(X)=\cup \overline{\{f(\C):f:\C \to X\}}$. Then $X\setminus Exc(X)$ is mordellic.
\end{enumerate}
\end{conjecture}

\textbf{Highlights in the history of the Green-Griffiths conjecture}

Here is a short (incomplete) list of results related to the Green-Griffiths conjecture.  

\begin{itemize}
\item In \cite{mcquillan2} McQuillen gave a positive answer to the conjecture for general surfaces if the second Segre class $c_1^2-c_2>0$ is positive. 
 \item In the seminal paper \cite{dem} Demailly---using ideas of Green, Griffiths and Bloch---works out a strategy for projective hypersurfaces. 
 \item In \cite{siu1,siu2} Siu gives positive answer for hypersurfaces of high degree, without effective lower bound for the degree.
 \item In \cite{dmr} Diverio, Merker and Rousseau give effective lower bound, proving that for a generic projective hypersurface of dimension $n$ and degree $>2^{n^5}$  the Green-Griffiths conjecture holds.
\item Recently, Merker (\cite{merker3}) has proved the existence of global jet differentials of high order for projective hyperpersurfaces in the optimal degree. Demailly in \cite{dem3} has proved the existence of global jet differentials (of possibly high order) for compact manifolds in general. 
\end{itemize}

\subsection{The strategy of the proof \cite{dem,dmr,dr,gg,siu2}}

Proving the algebraic degeneracy of holomorphic curves on $X$ means finding a non zero polynomial $P$ on $X$ such that all entire curves $f:\C \to X$ satisfy $P(f(\C))=0$. The first step of all known methods of proof establishes the existence of certain algebraic differential equations $P(f,f',\ldots,f^{(k)})=0$ of some order $k$, and then the second step is to find enough such equations so that they cut out a proper algebraic locus $Y\subset X$. 
Differential equations correspond to polynomial functions on the jet-bundles of holomorphic curves over $X$. These polynomial functions are called jet differentials, and their use can be traced back to the work of Bloch \cite{bloch}, Cartan \cite{cartan}, Ahlfors \cite{ahlfors}, Green and Griffiths \cite{gg}, Siu \cite{siu2}. Their ideas were extended in the seminal paper of Demailly \cite{dem}, by Diverio, Merker and Rousseau \cite{dmr} and by the author \cite{berczi}.

For more details on the history of this approach see \cite{dr,dem}.

Let 
\[f:\C \to X,\ \ t\to f(t)=(f_1(t),f_2(t), \ldots ,f_n(t))\]
be a curve written in some local holomorphic coordinates
$(z_1,\ldots ,z_n)$ on $X$. Let $J_kX$ be the $k$-jet bundle over $X$ of holomorphic curves, whose fibre at $x \in X$ is 
\[(J_kX)_x=\{\hat{f}_{[k]}: f:(\C,0) \to (X,x) \} \to X\]
sending $f_{[k]}$ to $f(0)$. This fibre is canonically isomorphic to $J_k(1,n)$.

The group of reparametrizations $\mathbf{G}_k=J_k^{\reg}(1,1)$
acts fibrewise on $J_kX$. The fibres of $J_kX$ can be identified with $J_k(1,n)$, and the action is linearised as in \eqref{group} before. Note that $\mathbf{G}_k=\C^* \ltimes U_d$ is a $\C^*$-extension of its maximal unipotent subgroup,
and for $\lambda \in \C^*$
\[(\lambda \cdot f)(t)=f(\lambda \cdot t),\text{ so } \lambda \cdot (f',f'',\ldots ,f^{(k)})=(\lambda f',\lambda^2 f'',\ldots ,\lambda^k f^{(k)}).\] 
Algebraic differential operators correspond to polynomial functions on $J_kX$, and we call these polynomial functions jet differentials, they have the form 
\[Q(f',f'',\ldots ,f^{(k)})=\sum_{\alpha_i \in \mathbb{N}^n}a_{\alpha_1,\alpha_2,\ldots \alpha_k}(f(t))(f'(t)^{\alpha_1}f''(t)^{\alpha_2}\cdots f^{(k)}(t)^{\alpha_k}),\]
where $a_{\alpha_1,\alpha_2,\ldots \alpha_k}(z)$ are holomorphic coefficients on $X$ and $t \to z = f(t)$  is a curve.

$Q$ is homogeneous of weighted degree $m$ under the $\C^*$ action if and only if 
\[Q(\lambda f',\lambda^2 f'',\ldots ,\lambda^k f^{(k)})=\lambda^m Q(f',f'',\ldots ,f^{(k)}).\]

\begin{defi}
\begin{itemize}
\item (\textbf{Green-Griffiths \cite{gg}}) Let $E_{k,m}^{GG}$ denote the sheaf of algebraic differential operators of order $k$ and weighted degree $m$.
 \item (\textbf{Demailly} \cite{dem}) The bundle of invariant jet differentials of order k and weighted
degree m is the subbundle $E_{k,m} \subset E_{k,m}^{GG}$, whose elements are
invariant under arbitrary changes of parametrization, i.e. for $\phi \in \mathbf{G}_k$
\[Q((f\circ \phi)',(f \circ \phi)'',\ldots ,(f \circ \phi)^{(k)})=\phi'(0)^mQ(f',f'',\ldots, f^{(k)}).\]
\end{itemize}
\end{defi}

We want to apply the general pinciple that for a $G$-space $X$ the ring of invariant functions on $X$ can be identified with polynomial functions on the quotient $X/G$. Roughly speaking we want 
\[\oplus_m (E_{k,m})_x =\oplus_m (E_{k,m}^{GG})^{\mathbb{U}}_x=\mathcal{O}((J_kX)_x)^{\mathbf{G}_k}=\mathcal{O}(J_k(1,n)/\mathbf{G}_k)\] 

Applying Theorem \ref{bercziszenes} fibrewise we get

\begin{prop}
\begin{enumerate}
\item The quotient $J_kX/\mathbf{G}_k$ has the structure of a locally trivial bundle over X, and there is a holomorphic embedding
\[\phi^{\PP}:J_kX/\mathbf{G}_k \hookrightarrow \PP(\wedge^k(\cotx \oplus \Sym^2(\cotx)
\oplus \ldots \oplus \Sym^k(\cotx)).\] 
 The fibrewise closure of the image $\calx_k=\overline{\mathrm{im}\phi^{\PP}}$
is a relative compactification of $J_k(\cotx)/\mathbf{G}_k$ over $X$.
 \item We have
\[(\pi_k)_*\calo_{\calx_k}(m)=\calo(E_{k,m{k+1 \choose 2}})\]
where $\pi_k:\PP(\wedge^k(\cotx \oplus \Sym^2(\cotx) \oplus \ldots
\oplus \Sym^k(\cotx)))\to X$ is the projection.
\end{enumerate}
\end{prop}

The strategy to solve the Green-Griffiths conjecture is based on the following

\begin{thm}[\textbf{Fundamental vanishing theorem
,Green-Griffiths '78, Demailly '95, Siu '96}]
Let $P\in H^0(X,E_{k,m} \otimes \mathcal{O}(-A))$ be a global algebraic
differential operator whose coefficients vanish on some
ample divisor A. Then for any $f:\C \to X$, $P(f_{[k]}(\C))\equiv 0$. (Note that $f_{[k]}(\C)\subset J_kX$.)
\end{thm}

\begin{cor} \begin{enumerate}
\item Let $\sigma$ be a nonzero element of 
\[H^0(\calx_{k},\calo_{\calx_{k}}(m) \otimes \pi^*\calo(-A))\simeq
H^0(X,E_{k,m{k+1 \choose 2}} \otimes \calo(-A)).\] 
Then $f_{[k]}(\C)\subset Z_\sigma$, where $Z_\sigma \subset \calx_d$ is the zero divisor of $\sigma$.  
 \item If $\{\sigma_j\}$ is a basis of global sections then the image $f(\C)$ lies in $Y=\pi_k(\bigcap Z_{P_j})$, hence the Green-Griffiths conjecture holds if there are
enough independent differential equations so that
$Y=\pi_k(\bigcap(Z_{P_j}))\varsubsetneq X.$
\end{enumerate}
\end{cor}

It is crucial to control in a more precise way the
order of vanishing of these differential operators along the ample divisor. Thus, we need here a slightly different theorem.

\begin{thm}[\textbf{\cite{dmr}}]
Assume that $n=k$, and there exist a $\delta=\delta(n) >0$ and $D=D(n,\delta)$ such that
\[H^0(\calx_{n},\calo_{\calx_{n}}(m) \otimes \pi^*K_X^{-\delta m})\simeq
H^0(X,E_{n,m{n+1 \choose 2}}\cotx \otimes K_X^{-\delta m})\neq 0 \]
 whenever $\deg(X)>D(n,\delta)$ provided that $m>m_{D,\delta,n}$ is large enough.  
Then the Green-Griffiths conjecture holds for 
\[\deg(X) \ge \max(D(n,\delta), \frac{n^2+2n}{\delta}+n+2).\]
\end{thm}

The goal is therefore to find a global section of $\calo_{\calx_{n}}(m) \otimes \pi^*K_X^{-\delta m}$ keeping $D(n,\delta)$ small. Following \cite{dmr}, we use the algebraic Morse inequalities of Trapani.  

\begin{thm}\cite{trap}  Let $L \to X$ be a holomorphic line bundle given as 
\[L=F \otimes G^{-1} \text{ with } F,G \text{ nef bundles.}\]
 Then for any nonnegative integer $q$ we have 
\[\sum_{j=0}^q(-1)^{q-j}h^j(X,L^{\otimes m}\otimes E)\le
r\frac{m^n}{n!}\sum_{j=0}^q(-1)^{q-j}{n \choose j}F^{n-j}\cdot
G^j+o(m^n).\] 
 Applying this with $q=1$ we get 
\begin{equation}\label{qequalsone}
F^n-nF^{n-1}G>0 \Rightarrow H^0(L^{\otimes m})\neq 0 \text{ for } m\gg 0.
\end{equation}
\end{thm}

In \cite{berczi} we prove that $F$ and $G$ are nef bundles in the following equality. 
\[\nonumber \underbrace{\calo_{\calx_{n}}(1)\otimes \pi^*K_X^{-\delta{n+1 \choose 2}}}_{L}=\underbrace{(\calo_{\calx_{n}}(1) \otimes
\pi^*\calo_X(2n^2))}_{F} 
\otimes \underbrace{(\pi^*\calo_X(2n^2)\otimes \pi^*K_X^{\delta{n+1 \choose 2}})^{-1}}_{G}.\]

Introduce the following notations:
\[h=c_1(\calo_X(1));\ c_1(K_X)=-c_1(X)=(d-n-2)h;\  
\calo_{\calx_n}(1)=\det \tau\]
where $\tau \to \calx_n$ is the tautological $n$-bundle.
Now $\dim(\calx_n)=n^2$, and according to \eqref{qequalsone} we want to prove the positivity of the following intersection number on $\calx_n$.
\[
I=\int_{\calx_n} (c_1(\det \tau)+2n^2\pi^*h)^{n^2}-
n^2(c_1(\det \tau)+2n^2\pi^*h)^{n^2-1}(2n^2\pi^*h+\delta {n+1 \choose 2} (d-n-2)h).\]

We apply localisation using the double fibration model \eqref{model} on the fibres of $\calx_n$. We need a stronger version of the vanishing property of the iterated residue to guarantee that only one fixed point's contribution is non-zero. After going through these technical difficulties in \cite{berczi} we arrive at 
\begin{prop}[Residue formula for the intersection number, \cite{berczi}]\label{residue}
\begin{multline}\nonumber
I=\int_X \mathrm{Res}_{\mathbf{z}=\infty} 
\frac{\prod_{i<j}(z_i-z_j)\,\mathcal{Q}_d(z_1\ldots z_n)R(\mathbf{z},h,d,\delta)}
{\prod_{1\le i+j \le l\le n}(z_i+z_j-z_l)(z_1\ldots z_n)^n}\cdot
\\
\cdot \prod_{l=1}^n\left(1+\frac{dh}{z_l}\right)\prod_{l=1}^n\left(1-\frac{h}{z_l}+\frac{h^2}{z_l^2}-\ldots \right)^{n+2}
\end{multline}
where
 \begin{multline}\nonumber
R(\mathbf{z},h,d,\delta)=(-z_1-\ldots -z_n+2n^2h)^{n^2}-\\
-n^2(-z_1-\ldots -z_n+2n^2h)^{n^2-1}(2n^2h+\delta {n+1 \choose 2}(d-n-2)h).
\end{multline}
\end{prop}

\noindent \textbf{Analysis of the formula}
 
\begin{itemize}
\item The iterated residue is the coefficient of $\frac{1}{z_1\ldots z_n}$, and has the form $h^np(d,n,\delta)$.
 \item Integration on $X$ is the substitution $h^n=d$, so the result is $dp(d,n,\delta)$.
 \item $I=p(n,d,\delta)=a_n(n,\delta)d^n+\ldots +a_0(n,\delta)$ is a degree-$n$ polynomial in $d=\deg(X)$, with polynomial coefficients in $n,\delta$.  
 \item The leading coefficient is 
\[a_n(n,\delta)=\left(1-n^2{n+1 \choose 2}\delta\right)\Theta(n),\]
where 
\[\Theta(n)=\text{constant term of }\frac{\displaystyle
 Q(\mathbf{z})\,\prod_{i<j}(z_i-z_j) (z_1+\ldots +z_n)^{n^2}}{\displaystyle
\prod_{i+j\le l \le n}
(z_i+z_j-z_l)\displaystyle (z_1\ldots z_n)^n }\]

\end{itemize}

In \cite{berczi} we prove that $\Theta(n)>0$ is positive, which implies  
\begin{cor}
For $\delta < \frac{2}{n^3(n+1)}$ the leading coefficient $a_n(n.\delta)$ of $I=p(n,d,\delta)$ is positive and therefore $I>0$ for $d\gg 0$ large enough. 
\end{cor}

The backgroud and experimental evidences of the following conjecture is explained in \cite{berczi}. It says that quotients of "neighbouring" coefficients of the Thom polynomial is polynomial.  
\begin{conjecture}[\cite{berczi}]\label{conj}
Define the Thom generating function $\mathrm{Tp}_k$ as 
\[\mathrm{Tp}_k(z_1,\ldots, z_k)=\frac{\prod_{m<l}(z_m-z_l)\,Q_k(z_1\ldots  z_k)}
{\prod_{m+r \le l\le k}(z_m+z_r-z_l)}.\]
Then 
\[\frac{\mathrm{coeff}_{z_1^{i_1} \ldots z_k^{i_k}}Tp_k}{\mathrm{coeff}_{z_1^{i_1} \cdots z_l^{i_l+1}\cdots  z_m^{i_m-1} \cdots z_k^{i_k}}Tp_k}<k^2\]
\end{conjecture}
Further analysis of the iterated residue formula in Proposition \ref{residue} leads us to 

\begin{thm}[\cite{berczi}]
Conjecture \ref{conj} and Conjecture \ref{rimanyiconj} for Thom polynomials of $A_n$ singularities implies the Green-Griffiths conjecture for generic smooth projective hypersurfaces $X\subset \P^{n+1}$ with $d=\deg(X)>n^6$.
\end{thm}

\section{Euler characteristic of the Demailly bundle}

The given iterated residue formula is suitable to compute other intersection numbers as well.  
The Euler characteristic of the Demailly bundle is defined as 
\[\chi(X,E_{k,m}T_X^*)=\sum_{i=0}^n (-1)^i \dim H^i(X,E_{k,m}T_X^*).\]
It is  well-known (see \cite{hirzebruch}) that 
\[\chi(X,E_{k,m})=\int_X [\mathrm{Ch}(E_{k,m})\cdot \mathrm{Td}(T_X)]_n\]
where $\mathrm{Ch}(\calo_{X_n}(1))$ is the Chern character and $\mathrm{Td}(T_X)=1+\frac{1}{2}c_1+\frac{1}{12}(c_1^2+c_2)+\ldots$ is the Todd-class.

\begin{thm}[\textbf{Iterated residue formula for the Euler-characteristics}]
\begin{multline*}
\chi(X,\pi_* \calo_{X_n}(m))=\\
\int_X \mathrm{Res}_{\mathbf{z}=\infty} 
\frac{\prod_{i<j}(z_i-z_j)\,\mathcal{Q}_n(z_1\ldots z_n)\mathrm{Ch}(\calo_{X_n}(m))\mathrm{Td}(T_X)}
{\prod_{1\le i+j \le l \le n}(z_i+z_j-z_l)(z_1\ldots z_n)^n}\cdot
\\
\cdot \prod_{l=1}^n\left(1+\frac{dh}{z_l}\right)\prod_{l=1}^n\left(1-\frac{h}{z_l}+\frac{h^2}{z_l^2}-\ldots \right)^{n+2}
\end{multline*}

where 
\[\mathrm{Ch}(\calo_{X_n}(1))=e^{m(z_1+\ldots +z_n)},\  \mathrm{Td}(T_X)=1+\frac{1}{2}c_1+\frac{1}{12}(c_1^2+c_2)+\ldots\]
\end{thm}

\section{Curviliear Hilbert schemes}

The goal of this section is to give a general framework for our localisation arguments. If $G_k=\jetreg 11$ denotes the group of $k$-jets of reparametrisation germs of $\C$ and $\jetreg 1n$ the $k$-jets of germs of curves $f:(\C,0) \to (\C^n,0)$, then the quotient 
$\jetreg 1n/\jetreg 11$ plays an important role in our applications, namely:

\begin{enumerate}
\item $\Sigma_k$ fibres over $\jetreg 1n/\jetreg 11$ with linear fibres. The Thom polynomials of Morin singularities are certain equivariant intersection numbers on $\Sigma_k$. 
\item $\jetreg 1n/\jetreg 11$ is isomorphic to the fibres of the Demailly jet bundle $E_{k}$ over a smooth manifold of dimension $n$. The positivity of the Demailly intersection number implies the Green-Griffiths conjecture.  
\end{enumerate}

In both applications we compute certain (equivariant) intersection numbers on the closure $\overline{\phi^{\mathrm{Grass}}(\jetreg 1n/ \jetreg 11)}\subset \grass_k(\oplus_{i=1}^k \Sym^i\C^n)$. 

In this subsection we identify the closure $\overline{\jetreg 1n/\jetreg 11}$ 
of $\jetreg 1n/\jetreg 11$ embedded in $\grass_k(\oplus_{i=1}^k \Sym^i\C^n)$  with the curvilinear component of the $k$-point punctual Hilbert scheme on $\C^n$; this geometric component of the punctual Hilbert scheme on $\C^n$ is thus the compactification of a non-reductive quotient.

Hilbert schemes of points on surfaces form a central object of geometry and representation theory and have a huge literature (see for example \cite{nakajima,bertin}). Recently many interesting connections between Hilbert schemes of points on  planar curve singularities and the topology of their links have been discovered \cite{shende,oblomkovshende,ors,maulik}.  However, much less is known about Hilbert schemes or punctual Hilbert schemes on higher dimensional manifolds. 

Let $(\C^n)_0^{[k]}$ be the punctual Hilbert scheme of $k$ points on $\C^n$, that is, the set of zero dimensional subschemes of $\C^n$ of length $k$ supported at the origin. The geometry of $(\C^n)_0^{[k]}$ is not knows; it is highly singular and irreducible. There is, however, an important distinguished irreducible component of $(\C^n)_0^{[k]}$, namely the punctual curvilinear Hilbert scheme, defined as follows 
\begin{defi}
Let $\mathfrak{m}=(x_1,\ldots, x_n)\subset \calo_{\C^n,0}$ denote the maximal ideal at the origin. The punctual curvilinear Hilbert scheme is defined as the closure of the set of curvilinear subschemes; that is, those which vanish on a curve up to order $n$:
\[\mathcal{C}_n^{[k]}=\{I \subset \mathfrak{m}:\mathfrak{m}/I \simeq t\C[t]/t^{k+1}\}.\]
\end{defi}

Note that 
\[\symk:=\mathfrak{m}/\mathfrak{m}^{k+1}=\oplus_{i=1}^k
\Sym^i\C^n\]
is the set of function-germs of degree $\le n$, and the punctual Hilbert scheme naturally sits in the Grassmannian
\[(\C^n)_0^{[k]} \subset \grass(k,\symk).\]
Looking at our embedding $\phi^\grass$ it is not hard to check (see \cite{bk2}) that 
\begin{prop} We have
\[\mathcal{C}_n^{[k]}=\overline{\phi^\grass(\jetreg 1n/\jetreg 11)}.\]
\end{prop}

This means that $\mathcal{C}_n^{[k]}$ can be described as a projective completion of a non-reductive quotient. When $n=2$ the punctual curvilinear component $\mathcal{C}_2^{[k]}$ is dense in $(\C^2)_0^{[k]}$, and therefore 
\begin{cor}We have 
\[(\C^2)_0^{[k]}=\overline{\phi^\grass(\jetreg 12/\jetreg 11)}.\]
\end{cor}

Our localisation method developed on $\overline{\jetreg 1n/\jetreg 11}$ is therefore allows us to  compute intersection numbers on the punctual curvilinear Hilbert scheme $\mathcal{C}_n^{[k]}$ for $k\le n$. A more detailed study of non-reductive quotients allows us to improve this technique, the details with more applications will be published later.

\end{document}